\documentclass{amsart}[12pt]
\usepackage{amsmath}
\usepackage{amstext}
\usepackage{amssymb}
\usepackage{amsthm}
\usepackage{amscd}

\newtheorem{anyprop}{Anyprop}[section]

\newtheorem{theorem}[anyprop]{Theorem}
\newtheorem{lemma}[anyprop]{Lemma}
\newtheorem{proposition}[anyprop]{Proposition}
\newtheorem{corollary}[anyprop]{Corollary}

\theoremstyle{definition}

\newtheorem{question}[anyprop]{Question}

\newtheorem{remark}[anyprop]{Remark}

\newcommand{\NN}{\mathbb{N}}
\newcommand{\ZZ}{\mathbb{Z}}
\newcommand{\QQ}{\mathbb{Q}}
\newcommand{\RR}{\mathbb{R}}
\newcommand{\CC}{\mathbb{C}}

\newcommand{\PP}{\mathbb{P}}

\newcommand  {\shF}     {\mathcal{F}}

\newcommand  {\shM}     {\mathcal{M}}

\newcommand  {\shL}     {\mathcal{L}}

\newcommand  {\shS}     {\mathcal{S}}
\newcommand  {\shT}     {\mathcal{T}}

\newcommand  {\foa}     {\mathfrak{a}}

\newcommand  {\fom}     {\mathfrak{m}}

\newcommand  {\dual}    {\vee}

\newcommand  {\Ext}     {\operatorname{Ext}}

\newcommand  {\length}  {\operatorname{length}}

\newcommand  {\lra}     {\longrightarrow}

\newcommand  {\modu}     {\operatorname{mod}}

\renewcommand{\O}       {\mathcal{O}}

\newcommand  {\Proj}    {\operatorname{Proj}}

\newcommand  {\ra}      {\rightarrow}

\newcommand  {\rk}    {\operatorname{rk}}

\newcommand  {\Spec}    {\operatorname{Spec}}

\newcommand  {\Syz}     {\operatorname{Syz}}

\newcommand{\comdots}{ , \ldots , }
\newcommand{\komdots}{ , \ldots , }
\newcommand{\plusdots}{ + \ldots + }

\theoremstyle{plain}

\theoremstyle{definition}

\theoremstyle{remark}

\newcommand{\ulim}{ {\rm ulim \,}}
\newcommand{\sheT}{\shT}
\newcommand{\ringR}{R}

\newcommand{\fieldK}{\mathbb{K}}
\newcommand{\leftc}{\left(}
\newcommand{\rightc}{\right)}
\newcommand{\modeq}{\equiv}

\begin{document}

\title[On the arithmetic of tight closure]{On the arithmetic of tight closure}
\author{Holger Brenner}
\author{Mordechai Katzman}
\address{Department of Pure Mathematics,
University of Sheffield, Hicks Building, Sheffield S3 7RH, United Kingdom\\
{\it Fax number}: 0044-114-222-3769}
\email{H.Brenner@sheffield.ac.uk, M.Katzman@sheffield.ac.uk}


\date{\today}


\begin{abstract}
We provide a negative answer to an old question in tight closure
theory by showing that the containment $x^3y^3 \in (x^4,y^4,z^4)^*$
in $\fieldK[x,y,z]/(x^7+y^7-z^7)$ holds for infinitely many but not
for almost all prime characteristics of the field $\fieldK$. This
proves that tight closure exhibits a strong dependence on the
arithmetic of the prime characteristic. The ideal $(x,y,z) \subset
\fieldK[x,y,z,u,v,w]/(x^7+y^7-z^7, ux^4+vy^4+wz^4+x^3y^3)$
has then the property that the
cohomological dimension fluctuates arithmetically between $0$ and
$1$.
\end{abstract}

\maketitle


\setcounter{section}{-1}
\section{Introduction}


\noindent This paper deals with a question regarding tight closure
in characteristic zero which we now review. Let $R$ be a commutative
ring of prime characteristic $p$ and let $I\subseteq R$ be an ideal.
Recall that for $e \geq 0$, the $e$-th Frobenius power of $I$,
denoted $I^{[p^e]}$, is the ideal of $R$ generated by all $p^{e}$-th
powers of elements in $I$. We say that $f\in I^*$, the tight closure
of $I$, if there exists a $c$ not in any minimal prime of $R$ with
the property that $c f^{p^e} \in I^{[p^e]}$ for all large $e\geq 0$.
This notion, due to M.~Hochster and C.~Huneke, is now an important
tool in commutative algebra and algebraic geometry, particularly
since it gives a systematic framework for reduction to positive
characteristic. We refer the reader to \cite{hochsterhunekebriancon}
for the basic properties of tight closure in characteristic $p$.

\smallskip\noindent
How does the containment $f \in I^*$ depend on the prime
characteristic? To make sense of this question suppose that
$R_{\ZZ}$ is a finitely generated ring extension of $\ZZ$ and that
$I \subseteq R_{\ZZ}$ is an ideal, $f \in R_{\ZZ}$. Then we may
consider for every prime number $p$ the specialization $R_{\ZZ/(p)}=
R_{\ZZ} \otimes_{\ZZ} \ZZ/(p)$ of characteristic $p$ together with
the extended ideal $I_{p} \subseteq R_{\ZZ/(p)}$, and one may ask
whether $f_p \in I_{p}^*$ holds or not. We refer to this question
about the dependence on the prime numbers as the `arithmetic of
tight closure'.

\smallskip\noindent
Many properties in commutative algebra exhibit an arithmetically
nice behaviour: for example, $R_\QQ$ is smooth (normal,
Cohen-Macaulay, Gorenstein) if and only if $R_{\ZZ/(p)}$ is smooth
(normal, Cohen-Macaulay, Gorenstein) for almost all prime numbers
(i.e., for all except for at most finitely many). In a similar way
we have for an ideal $I \subseteq R_{\ZZ} $ that $I_{\QQ}=I R_{\QQ}$
is a parameter ideal or a primary ideal if and only if this is true
for almost all specializations $I_p$. Furthermore, $f \in I$ if and
only if $f_p \in I_p$ holds for almost all prime characteristics:
see \cite[Chapter 2.1]{hochsterhuneketightzero} and appendix 1 in
\cite{hunekeapplication} for this kind of results.

\smallskip\noindent
When $R$ is a finitely generated $\mathbb{\QQ}$-algebra, Hochster
and Huneke define the tight closure of an ideal $I\subseteq R$, in
the same spirit as the examples above, with the help of a
$\ZZ$-algebra $R_{\ZZ}$ where $R= R_{\ZZ} \otimes_{\ZZ} \QQ$, as the
set of all $f \in R$ for which $f _p \in (I_p)^*$ holds for almost
all $p$. This definition is independent of the chosen model
$R_{\ZZ}$. The reader should consult \cite{hochsterhuneketightzero}
for properties of tight closure in characteristic zero. This
definition works well, because the most important features from
tight closure theory in positive characteristic, like $F$-regularity
of regular rings, colon capturing, Brian\c{c}on-Skoda theorems,
persistence, behave well arithmetically, so that these properties
pass over to the characteristic zero situation with full force.

\smallskip\noindent
M.~Hochster and C.~Huneke (see appendix 1 in
\cite{hunekeapplication} or Question 11 in the appendix of
\cite{hochsterhuneketightzero} or Question 13 in
\cite{hochstertightsolid}) and the second author (see \S 4 in
\cite{katzmanfinitecriteria}) raise the following natural question:
if $R$ is a finitely generated $\mathbb{Z}$-algebra of
characteristic zero and $I\subseteq R$ is an ideal which is tightly
closed, i.e. $I^*=I$ in $R_\QQ$, must one have $(I_p)^*=I_p$ for
almost all primes $p$? Or, using the terminology of
\cite{katzmanfinitecriteria}, must tightly closed ideals be
fiberwise tightly closed?

\smallskip\noindent
As often in tight closure theory, the situation for parameter ideals
is better understood than the general case, but even for parameter
ideals a complete answer is not known. There are however results due
to N.~Hara and K.~Smith (see \cite{hararationalfrobenius},
\cite{haratestmultiplier}, \cite[Theorem 6.1]{hunekeparameter},
\cite{smithtestmultiplier}, \cite[Theorem 2.10, Open Problem
2.24]{smithvanishing}) which imply that for a normal standard-graded
Cohen-Macaulay domain with an isolated singularity and for a normal
Gorenstein algebra of finite type over a field the answer is
affirmative.

\smallskip\noindent
The main theorem in this paper (Theorem \ref{mainTheorem}) provides,
however, a negative answer to this question by showing that for the
homogeneous primary ideal $I=(x^4,y^4,z^4)$ in
$\mathbb{Z}[x,y,z]/(x^7+y^7-z^7)$ one has $x^3 y^3 \in (I_p) ^*$ for
$p \modeq 3 \modu 7$ but $x^3 y^3 \notin (I_p)^*$ for $p \modeq 2
\modu 7$.

\smallskip\noindent
Our example has also interesting implications for the dependence of
the cohomological dimension on the characteristics of ground fields.
The ideal $\foa =(x,y,z)$ inside the forcing algebra
$A=\fieldK[x,y,z,u,v,w]/(x^7+y^7-z^7, ux^4+vy^4+wz^4+x^3y^3)$ is
such that the open subset $D(\foa ) \subset \Spec A$ is affine for
infinitely many but not for almost all prime reductions. This means
that its cohomological dimension fluctuates arithmetically between
$0$ and $1$, see \ref{solidclosure} for this relation via solid
closure and \ref{cohodimprojective} for an interpretation in terms
of projective varieties.

\smallskip\noindent
Moreover, our example has also consequences for the study of
Hilbert-Kunz multiplicities which we discuss in \ref{hilbertkunz}
and for the non-standard tight closure of H.~Schoutens (see
\ref{nonstandard}).


\smallskip\noindent
During the preparation of this paper we used the computer algebra
systems Cocoa and Macaulay 2 (\cite{cocoaSystem}, \cite{M2}). We
thank A. Kaid and R. Y. Sharp for useful communications.

\section{Reduction to Frobenius powers}

\noindent In this section we show where to look for candidates
$(R,I,f)$ with the property that $f_p \in I_p^*$ holds for
infinitely many but not for almost all prime numbers $p$ . This
approach rests on the geometric interpretation of tight closure in
terms of bundles, which we now recall briefly. Let $\ringR$ denote a
geometrically normal two-dimensional standard-gra\-ded domain over a
field $\fieldK$. A set of homogeneous generators $f_1 \komdots f_n
\in \ringR$ of degrees $d_1 \komdots d_n$ of an $\ringR_+$-primary
ideal give rise to the short exact sequence of locally free sheaves
on the smooth projective curve $C= \Proj \ringR$,
$$0 \lra \Syz(f_1 \komdots f_n)(m) \lra \bigoplus_{i=1}^n \O_C(m-d_i)
\stackrel{f_1 \komdots f_n}{\lra} \O_C(m) \lra 0 \, .$$ A
homogeneous element $f \in \ringR$ of degree $m$ defines via the
connecting homomorphism a cohomology class $\delta(f) \in
H^1(C,\Syz(f_1 \komdots f_n)(m))$ in this syzygy sheaf. It was shown
in \cite{brennertightproj}, \cite{brennerslope} how this cohomology
class is related to the question as to whether $f$ belongs to the
tight closure (in positive characteristic) of the ideal $(f_1
\komdots f_n)$ or not. The cohomology class $c \in H^1(C, \shS)=
\Ext^1(\O_C, \shS)$ corresponds to an extension $0 \ra \shS \ra
\shS' \ra \O_C \ra 0$ and to a geometric torsor $\PP(\shS'^\dual)-
\PP(\shS^\dual)$. Now $f \in (f_1 \komdots f_n)^*$ if and only if
the torsor defined by $\delta(f)$ is not an affine scheme.

\smallskip\noindent
If the syzygy bundle is strongly semistable in positive
characteristic $p$, then this approach gives a numerical criterion
for $(f_1 \komdots f_n)^*$, where the degree bound which separates
inclusion from exclusion is given by $(d_1 \plusdots d_n)/(n-1)$.
So, if we want to find an example where $f \in (f_1 \komdots f_n)^*$
holds for infinitely many prime numbers but not for almost all, we
have to look first for an example where for infinitely many prime
numbers the syzygy bundle is not strongly semistable (this is also
the reason why such an example cannot exist in the cone over an
elliptic curve). That $\shS$ is not strongly semistable means that
some Frobenius pull-back of it, say $\shT= F^{e*}(\shS)$, is not
semistable, and that means that there exists a subbundle $\shF
\subset \shT$ such that $\deg(\shF)/ \rk (\shF)  > \deg (\shT)/
\rk(\shT)$.

\smallskip\noindent
Examples of such syzygy bundles with the property that they are
semistable in characteristic zero but not strongly semistable for
infinitely many prime numbers were first given in
\cite{brennermiyaoka}, where it was shown that a question of Miyaoka
and Shepherd-Barron (\cite{miyaokachern},
\cite{shepherdbarronsemistability}) has a negative answer. The
following lemma gives another example of that kind.

\begin{lemma}
\label{notsemistable} Let $d \in \NN$ and let $p$ denote a prime
number; write $p= d \ell +r $, $0 < r <d$. Suppose that $ d/4 \leq r
< d/3$. Let $\fieldK$ denote a field of characteristic $p$ and let
$C=\Proj \fieldK[x,y,z]/(x^d+y^d-z^d)$ be the Fermat curve of degree
$d$. Then the first Frobenius pull-back of $\Syz(x^{4},y^{4},z^{4})$
on $C$ is not semistable.
\end{lemma}
\proof
We have $4p= 4d \ell +4r = d(4 \ell+1) +(4r - d)$; set
$t=4r-d$. We consider first in $\fieldK[x,y]$ the syzygies for
$$  x^{4p}=x^{d(4 \ell+1) + t},\, y^{4p}= y^{d(4\ell+1) +t},\,  (x^d+y^d)^{4\ell+1} \, .$$
We multiply the last term by the $2\ell+1$ monomials
$$x^ty^t(x^{d2\ell}y^{0}), \,x^ty^t(x^{d(2\ell-1)}y^{d})
\comdots x^ty^t(x^0y^{d2 \ell}) \, .$$ The resulting polynomials are
expressible modulo the first two terms as a $\fieldK$-linear
combination of the monomials $x^ty^t x^{di}y^{ d j}$, where $i+j = 6
\ell+1$ and $i,j \leq 4\ell$. Therefore $i=2\ell+1 \komdots 4\ell$
and there are only $2\ell$ of these. Hence there exists a global
non-trivial syzygy $(h_1,h_2,h_3)$ of these polynomials of total
degree $d(6\ell+1) +2t$. Therefore $(z^th_1, z^th_2,h_3)$ is a
global non-trivial syzygy for $x^{4p},y^{4p},z^{4p}$, since $$ 0=
z^th_1x^{4p} +z^th_2y^{4p}+z^th_3(x^d+y^d)^{4 \ell +1}= z^th_1x^{4p}
+z^th_2y^{4p}+h_3 z^{4p} \, .$$ The total degree of this syzygy is
$d(6\ell+1)+3t$. The degree of the bundle
$$\Syz(x^{4p},y^{4p},z^{4p})(d(6\ell+1) +3t)$$
is however (up to the factor $\deg (\O(1))$)
$$2( d(6\ell+1)+3t) -3(d(4\ell+1 ) +t) = -d+3t = 12r -4d       \, ,$$
which is negative due to the assumption that $r < d/3$. But a bundle
of negative degree and with a non-trivial section is not semistable.
\qed

\medskip\noindent
The following proposition reduces under suitable conditions the
computation of tight closure to the computation of a certain
Frobenius power.

\begin{proposition}
\label{secondfrobenius} Let $\fieldK$ denote a field of positive
characteristic $p$ and let $\ringR$ denote a two-dimensional
geometrically normal standard-graded domain over $\fieldK$. Suppose
that $p \geq 2g+1$, where $g$ denotes the genus of the smooth
projective curve $C= \Proj \ringR$. Let $f_1,f_2,f_3$ denote
homogeneous elements in $\ringR$ which generate an
$\ringR_+$-primary ideal. Let $m \in \ZZ$ be such that the $e$-th
pull-back of the syzygy bundle $\Syz(f_1,f_2,f_3)(m)$ can be
incorporated in a short exact sequence on $C$,
$$ 0 \lra \shL \lra F^{e*} ( \Syz(f_1,f_2,f_3)(m))= \Syz(f_1^{q},f_2^{q},f_3^{q})(qm)
\lra \shM \lra 0 \, ,$$ where $q=p^{e}$ and $\shL$ is an invertible
sheaf of positive degree and $\shM$ is an invertible sheaf of
negative degree. Let $f$ denote a homogeneous element of degree $m$.
Then $f \in (f_1,f_2,f_3)^*$ if and only if $f^{pq} \in (
f_1^{pq},f_2^{pq},f_3^{pq})$.
\end{proposition}
\proof The implication from right to left is clear. For the other
direction we may assume that $\fieldK$ is algebraically closed. We
will argue on the smooth projective plane curve $C=\Proj \ringR$ and
use the geometric interpretation of tight closure. We apply the
Frobenius to the given short exact sequence and obtain a new exact
sequence
$$0 \lra \shL^p \lra \Syz(f_1^{pq},f_2^{pq}, f_3^{pq})(pqm) \lra \shM^p \lra 0 \, .$$
The cohomology sequence is
$$ \lra H^1(C,\shL^p) \lra
H^1(C,\Syz(f_1^{pq},f_2^{pq}, f_3^{pq})(pqm)) \lra H^1(C,\shM^p)
\lra 0.$$ The genus of the curve $C$ is $g$ and the canonical sheaf
$\omega_C$ has degree $2g-2$. Hence for $p> 2g-2$ we have that $\deg
(\shL^{-p} \otimes \omega_C) <0$ and therefore $H^1(C, \shL^p)=0$ by
Serre duality. This gives an isomorphism
$$  H^1(C,\Syz(  f_1^{pq},f_2^{pq}, f_3^{pq})(pqm))\cong H^1(C,\shM^p) \, .$$
Suppose now that $f^{pq} \not\in (f_1^{pq}, f_2^{pq},f_3^{pq})$.
This means that the corresponding cohomology class $c=\delta
(f^{pq}) \in H^1(C,\Syz(f_1^{pq},f_2^{pq}, f_3^{pq})(pqm))$ is not
zero; let $c' \neq 0$ denote the corresponding class in $H^1(C,
\shM^p)$. To show that $f$ does not belong to the tight closure of
$(f_1,f_2,f_3)$ we show that the geometric torsor corresponding to
$c$ is an affine scheme \cite[Proposition 3.9]{brennertightproj},
and for that it is sufficient to show that the geometric torsor
corresponding to $c'$ is an affine scheme. The class $c'\in H^1(C,
\shM^p) \cong \Ext^1(\O_C, \shM^p)$ defines a non-trivial extension
$$ 0 \lra \shM^p \lra \sheT \lra \O_C \lra 0$$
with dual sequence
$$ 0 \lra \O_C \lra \sheT^\dual \lra \shM^{-p} \lra 0 \, .$$
Here $\shM^{-p}$ is ample, since its degree is positive, and
therefore by \cite[Proposition 2.2]{giesekerample} every quotient
bundle of $\sheT^\dual$ has positive degree. Since $\deg \sheT^\dual
= \deg \shM^{-p} \geq p > 2 \cdot g$, it follows by \cite[Lemma
2.2]{giesekerample} that $\sheT^\dual$ is an ample vector bundle
(one can also argue using \cite[Corollary 7.7]{hartshorneample}).
But then $C \cong \PP(\shM^{-p})\subset\PP(\sheT^\dual)$ is an ample
divisor and its complement is affine. \qed

\begin{remark}
\label{semistableremark}
The situation described in Proposition
\ref{secondfrobenius} occurs in particular for
$$2m= \deg(f_1)+\deg(f_2)+ \deg(f_3)$$
under the condition that the syzygy bundle is
not strongly semistable. For then some Frobenius pull-back
$\shT=\Syz(f_1^q,f_2^q,f_3^q)(qm)$ is not semistable, but its degree
is
$$q(2m-\deg(f_1)-\deg(f_2)-\deg(f_3))\deg (\O_C(1))=0 \, .$$ Then
there exists the maximal destabilizing invertible subsheaf $\shL
\subset \shT$ of positive degree, and the quotient sheaf is also
invertible of negative degree.
\end{remark}

\begin{corollary}
\label{secondfrobenius27}
Let $\fieldK$ denote a field of positive
characteristic $p \modeq 2 \modu 7$, $ p \neq 2, 23$, and let
$\ringR=\fieldK[x,y,z]/(x^7+y^7-z^7)$. Then $x^3y^3 \in
(x^4,y^4,z^4)^*$ if and only if $x^{3p^2}y^{3p^2} \in
(x^{4p^2},y^{4p^2},z^{4p^2})$.
\end{corollary}
\proof
In the notation of Lemma \ref{notsemistable} we have $p=d
\ell+2$ hence $r=2$ and clearly $7/4 \leq 2 < 7/3$. Hence the first
Frobenius pull-back of $\Syz(x^4,y^4,z^4)$ is not semistable.
Therefore via Remark \ref{semistableremark} we are in the situation
of Proposition \ref{secondfrobenius} with $e=1$ hence $pq=p^2$.
Since $g=15$, the condition on the prime number is $p \geq 31$, so
only $p= 2$ and $23$ are excluded. \qed

\begin{remark}
\label{2,23remark} The method of Proposition \ref{secondfrobenius}
works in principle also for small prime numbers $p$. We only need to
find a power $p^u \geq 2g+1$. If the $e$-th Frobenius pull-back is
not semistable, then we can conclude that $f \in (f_1,f_2,f_3)^*$ if
and only if $f^{p^{u}p^{e}} \in
(f_1^{p^{u}p^{e}},f_2^{p^{u}p^{e}},f_3^{p^{u}p^{e}})$. So in
Corollary \ref{secondfrobenius27} we take $u=2$ for $p=23$ and $u=5$
for $p=2$ to make things work also in these cases.
\end{remark}

\section{The case $p\modeq 2 \modu 7$}

\noindent In this section we want to show that $x^3y^3 \not\in
(x^4,y^4,z^4)^*$ in $\fieldK[x,y,z]/(x^7+y^7-z^7)$ if $\fieldK$ has
characteristic $p\modeq 2 \modu 7$. We will need the following
lemmata on matrices.

\begin{lemma}
\label{lemma1}
The $r\times s$ matrix $A$ with entries
$$\leftc {a\choose{b+i-j}}\rightc_{1\leq i\leq r, 1\leq j\leq s}$$
can be brought to the form
$$\leftc {a+j-1\choose{b+i-1}}\rightc_{1\leq i\leq r, 1\leq j\leq s}$$
by performing elementary column operations.
\end{lemma}
\begin{proof}
We proceed by induction on $s$. If $s=1$ there is nothing to show,
so assume $s>1$. Add the penultimate column of $A$ to the last
column; in the result, add the $(s-2)$-th column to the $(s-1)$-th,
and continue in this way until the first column has been added to
the second column. In this way one obtains the matrix
$$
\left(
\begin{array}{ccccc}
  {{a} \choose{b}} & {{a+1} \choose{b}} &{{a+1} \choose{b-1}} & \dots & {{a+1} \choose{b+2-s}} \\
  {{a} \choose{b+1}} & {{a+1} \choose{b+1}} &{{a+1} \choose{b}} & \dots & {{a+1} \choose{b+3-s}} \\
  \vdots &  \vdots & & & \vdots  \\
  {{a} \choose{b+r-1}} & {{a+1} \choose{b+r-1}}  & {{a+1} \choose{b+r-2}} & \dots & {{a+1} \choose{b+r-s+1}}
\end{array}\right).
$$
Now apply the induction hypothesis to the submatrix of this matrix consisting of all its
 columns except the first.
\end{proof}

\medskip\noindent
Using Lemma \ref{lemma1} one can obtain the following result due to
V.~van Zeipel (\cite{Zeipel}; the calculation is described in
\cite[Chapter XX]{Muir}.)
\begin{lemma}
\label{lemma2}
$$ \det \leftc {{a} \choose{b+i-j}} \rightc_{1\leq i\leq r \atop {1\leq j\leq r}}=
\prod_{t=0}^{r-1} \frac{{{a+r-1-t} \choose{b}}}{{{b+t} \choose{b}}}
\, .
$$
\end{lemma}


\medskip\noindent
We will use these lemmata in the proof of the following result.

\begin{lemma}
\label{lemmafrobenius1}
Let $\fieldK$ denote a field of positive
characteristic $p=7 \ell +2$. Then we have $x^{3p}y^{3p} \not\in
(x^{4p}, y^{4p}, (x^7+y^7)^{4\ell+1})$ in the polynomial ring
$\fieldK[x,y]$.
\end{lemma}
\begin{proof}
The case $p=2$ is checked immediately, so suppose that $\ell > 0$.
Since $3p= 21 \ell +6$ and $4p= 28 \ell +8$, we rewrite what we want
to show as
\begin{equation}\label{eq1}
x^{21 \ell +6} y^{21 \ell +6} \notin (x^{28 \ell+ 8}, y^{28 \ell
+8}, (x^7+y^7)^{4 \ell+ 1}) \,  .
\end{equation}
We endow $ \fieldK[x,y]$ with a $\mathbb{Z}/7\mathbb{Z} \oplus
\mathbb{Z}/7\mathbb{Z} \oplus \mathbb{Z}$ grading by assigning $x$
degree $(\overline{1},\overline{0},1)$ and $y$ degree
$(\overline{0},\overline{1},1)$. With this grading the left hand
side of (\ref{eq1}) is homogenous of degree
$(\overline{6},\overline{6},42 \ell+12)$ while $(x^7+y^7)^{4 \ell+
1}$ is homogeneous of degree $(\overline{0},\overline{0},28 \ell
+7)$, so the degree difference is $(\overline{6},\overline{6},14
\ell +5)$. Condition (\ref{eq1}) fails to hold if and only if there
exist $a_0, \dots, a_{2 \ell -1}  \in \fieldK$ such that
\begin{equation*}
x^{21 \ell +6} y^{ 21 \ell +6} \equiv \big(\sum_{i=0}^{2 \ell -1}
a_i x^{7i+6} y^{7(2\ell-1-i)+6} \big) \big(\sum_{j=0}^{4 \ell+1} {4
\ell+1\choose{j}} x^{7j} y^{7(4 \ell+1-j)} \big) \mod(x^{28 \ell
+8},y^{28 \ell+ 8})
\end{equation*}
and we assume that this is the case. Notice that the terms occurring
on the right hand side of this equation have the form
$x^{7i+6}y^{7(6 \ell -i)+6}$ for $0\leq i \leq 6 \ell $. Since
$$
\left\{
\begin{array}{lll}
7i+6& < &28 \ell +8\\
7(6 \ell -i)+6& < &28 \ell + 8
\end{array}
\right. \Leftrightarrow\quad 2 \ell  \leq i\leq 4 \ell \, ,
$$
we obtain $\modu (x^{28 \ell+8},y^{28 \ell+8})$
\begin{eqnarray*}
x^{21 \ell+6} y^{21 \ell +6} & \equiv & x^6 y^6 \sum_{i=2 \ell}^{4 \ell}
\left( \sum_{j=0}^{2\ell -1} a_j {4 \ell +1\choose{i-j}} \right) x^{7i}y^{7(6 \ell-i)}
\\
& \equiv & x^6 y^6 \sum_{i=1}^{2 \ell +1} \left( \sum_{j=1}^{2 \ell}
a_{j-1} {4 \ell +1 \choose{2 \ell +i-j}} \right) x^{7(2 \ell
+i-1)}y^{7(4 \ell -i+1)}
\end{eqnarray*}
and since no term in the last expression is divisible by $x^{28 \ell
+8}$ or by $y^{28 \ell +8}$, we deduce that
$$
x^{21 \ell+6} y^{21 \ell +6} = x^6 y^6 \sum_{i=1}^{2 \ell +1} \left(
\sum_{j=1}^{2 \ell} a_{j-1} {4 \ell +1 \choose{2 \ell +i-j}} \right)
x^{7(2 \ell +i-1)}y^{7(4 \ell -i+1)} \, .$$ We may cancel $x^6y^6$
from both sides of the equation and we write $X=x^7, Y=y^7$ to
obtain
\begin{equation*}
X^{3 \ell } Y^{3 \ell } = \sum_{i=1}^{2 \ell +1} \left(
\sum_{j=1}^{2 \ell} a_{j-1} {4 \ell +1 \choose{2 \ell +i-j}} \right)
X^{2 \ell +i-1}Y^{4 \ell -i+1} \, .
\end{equation*}
If we compare the coefficients of $X^{2\ell +i -1}Y^{4 \ell +1-i}$
for $1 \leq i \leq 2 \ell +1$ we obtain the conditions
$$\sum_{j=1}^{2 \ell} a_{j-1} {4 \ell + 1 \choose{2 \ell +i-j}}
= \delta_{i, \ell+1} \,\,\, \mathrm{\ for\ all\ } \, \,\, 1 \leq i
\leq 2\ell+1
$$ where $\delta_{i, \ell +1}$ is Kronecker's delta. If we define $M_1$
to be the $(2 \ell +1 ,2 \ell)$ matrix whose entries are
$${4 \ell + 1 \choose{2 \ell +i-j}}_{1\leq i\leq 2 \ell +1, 1\leq j\leq 2 \ell}$$
and if $\mathbf{e}_{\ell+1}$ is the $( \ell+1)$-th elementary column
vector of  size $2 \ell +1$ then we are now assuming that
$\mathbf{e}_{\ell +1}$ is in the span of the columns of $M_1$. We
have to show that this is not possible. Since $M_1$ has more rows
than columns, its rows are linearly dependent, i.e., there exists a
$\mathbf{\rho}=(\rho_1, \dots, \rho_{2\ell+1})\neq \mathbf{0}$ such
that $\mathbf{\rho} M_1=0$. It is now enough to show that we can
choose this $\mathbf{\rho}$ with $\rho_{\ell+1}\neq 0$, since for
such $\mathbf{\rho}$ we have $\mathbf{\rho}
\mathbf{e}_{\ell+1}=\rho_{\ell+1}\neq 0$ and so
$\mathbf{e}_{\ell+1}$ could not be in the span of the columns of
$M_1$. Assume by way of contradiction that we can find a non-zero
$\mathbf{\rho}$ as above with  $\rho_{\ell+1}=0$. This implies that
the rows of $M_1$ numbered $1,\dots, \ell,\ell+2, \dots, 2\ell+1$
are linearly dependent. Use Lemma \ref{lemma1} and apply elementary
column operations to $M_1$ to obtain the matrix
\begin{equation*}
M_2= \leftc {4\ell+j\choose{2\ell+i-1}}  \rightc_{1\leq i\leq
2\ell+1, \atop{ 1\leq j\leq 2\ell}} = \left(
\begin{array}{llll}
4\ell+1 \choose{2\ell} & 4\ell+2 \choose{2\ell} & \dots & 6\ell \choose{2\ell} \\
4\ell+1 \choose{2\ell+1} & 4\ell+2 \choose{2\ell+1} & \dots & 6\ell \choose{2\ell+1} \\
\vdots & \vdots & & \vdots \\
4\ell+1 \choose{4\ell} & 4\ell+2 \choose{4\ell} & \dots & 6\ell \choose{4\ell}
\end{array}
\right) .
\end{equation*}
Use the fact that $\displaystyle {a+1 \choose{b+1}} =
\frac{a+1}{b+1} {a\choose{b}}$ and multiply rows $1,2, \dots
2\ell+1$ by $2\ell, 2\ell+1, \dots, 4\ell$ and divide columns $1,2,
\dots 2\ell$ by $4\ell+1, 4\ell+2, \dots, 6\ell$  to obtain
\begin{equation*}
M_2= \Lambda\big( \frac{1}{2\ell}, \frac{1}{2\ell+1}, \dots ,
\frac{1}{4\ell}\big ) \left(
\begin{array}{llll}
4\ell \choose{2\ell-1} & 4\ell+1 \choose{2\ell-1} & \dots & 6\ell-1 \choose{2\ell-1} \\
4\ell \choose{2\ell} & 4\ell+1 \choose{2\ell} & \dots & 6\ell-1 \choose{2\ell} \\
\vdots & \vdots & & \vdots \\
4\ell \choose{4\ell-1} & 4\ell+1 \choose{4\ell-1} & \dots & 6\ell-1 \choose{4\ell-1}
\end{array}
\right) \Upsilon\big( {4\ell+1}, {4\ell+2}, \dots, {6\ell} \big)
\end{equation*}
where
$\Lambda\big( a_1, \dots , a_{2\ell+1} \big )$ is the $(2\ell+1) \times (2\ell+1)$ diagonal matrix with
$a_1, \dots , a_{2\ell+1}$ along its diagonal and
$\Upsilon\big( b_1, \dots , b_{2\ell} \big )$ is the $2\ell \times 2\ell$ diagonal matrix with
$b_1, \dots , b_{2\ell}$ along its diagonal.
We can repeat this process $2\ell$ times to obtain
\begin{equation}\label{eq2}
M_2= \Lambda \left(
\begin{array}{llll}
2\ell+1 \choose{0} & 2\ell+2 \choose{0} & \dots & 4\ell \choose{0} \\
2\ell+1 \choose{1} & 2\ell+2 \choose{1} & \dots & 4\ell \choose{1} \\
\vdots & \vdots & & \vdots \\
2\ell+1 \choose{2\ell} & 2\ell+2 \choose{2\ell} & \dots & 4\ell \choose{2\ell} \\
\end{array}
\right)
\Upsilon
\end{equation}
where
$$\Lambda=\left(
\prod_{t=0}^{2\ell-1}\Lambda\big( \frac{1}{2\ell-t},
\frac{1}{2\ell+1-t}, \dots , \frac{1}{4\ell-t}\big ) \right)$$
$$\Upsilon=\left(
\prod_{t=0}^{2\ell-1}\Upsilon\big({4\ell+1-t}, {4\ell+2-t}, \dots,
{6\ell-t} \big) \right)
$$
We notice that none of the entries in the diagonal matrices above is
$0$ or $1/0$ modulo $p$ and so, if we denote with $M_3$ the middle
matrix in Equation \ref{eq2}, and, if we write
$\mathbf{\rho}^\prime=\mathbf{\rho} \Lambda$ then $\mathbf{\rho}
M_2= 0 \Leftrightarrow \mathbf{\rho}^\prime M_3=0$, $\mathbf{\rho}=0
\Leftrightarrow \mathbf{\rho}^\prime=0$ and $\rho_{\ell+1}=0
\Leftrightarrow \rho_{\ell+1}^\prime=0$. It is now, therefore,
sufficient to show that the rows of $M_3$ numbered $1,\dots,
\ell,\ell+2, \dots, 2\ell+1$ are not linearly dependent. Use Lemma
\ref{lemma1} to perform the inverse elementary column operations on
$M_3$ to bring it to the form
\begin{equation*}
M_4= \leftc { 2\ell+1 \choose{i-j} } \rightc_{1\leq i\leq 2\ell+1,
\atop{ 1\leq j\leq 2\ell}} = \left(
\begin{array}{lllll}
2\ell+1 \choose{0} & 0 & 0 & \dots & 0 \\
2\ell+1 \choose{1} & 2\ell+1 \choose{0} & 0 & \dots & 0 \\
 & \ddots & & & \\
2\ell+1 \choose{2\ell-1} & 2\ell+1 \choose{2\ell-2} & {2\ell+1\choose{2\ell-3}} & \dots & {2\ell+1\choose{0}} \\
2\ell+1 \choose{2\ell} & 2\ell+1 \choose{2\ell-1} & {2\ell+1\choose{2\ell-2}} & \dots & {2\ell+1\choose{1}} \\
\end{array}
\right)
\end{equation*}
and we now need to show that the rows of $M_4$ numbered $1,\dots,
\ell,\ell+2, \dots, 2\ell+1$ are not linearly dependent. If we
delete the $(\ell+1)$-th row from $M_4$ and perform elementary row
operations consisting of adding multiples of rows $1,2, \dots, \ell$
to lower rows, we can bring the resulting matrix to the form
$I_\ell\oplus M_5$ where $I_\ell$ is a $\ell\times \ell$ identity
matrix and $M_5$ consists of the lower-rightmost block of size
$\ell\times \ell$ in $M_4$, i.e.,
$$
M_5= \left(
\begin{array}{llll}
2\ell+1 \choose{1} & 2\ell+1 \choose{0}  & \dots & 0 \\
 & \ddots & & \\
2\ell+1 \choose{\ell-1} & 2\ell+1 \choose{2\ell-2} & \dots & {2\ell+1\choose{0}} \\
2\ell+1 \choose{\ell} & 2\ell+1 \choose{\ell-1} & \dots & {2\ell+1\choose{1}} \\
\end{array}
\right) .
$$
The value of the determinant of $M_5$ can be computed using Lemma
\ref{lemma2}:
$$\det M_5= \prod_{t=0}^{\ell-1} \frac{\binom{2 \ell +1 +\ell -1 -t}{1}}{\binom{1+t}{1}}=
\prod_{t=0}^{\ell-1} \frac{3\ell-t}{1+t} $$ which is a unit modulo
$p$. Hence the rows of $M_5$ are linearly independent, and we
conclude that $\mathbf{e}_{\ell+1}$ in not in the span of the
columns of $M_2$.
\end{proof}

\begin{proposition}
\label{theorem 2 mod 7}
If  $p \modeq 2 \modu 7$, then $x^3 y^3 \notin (x^4,y^4,z^4)^*$ in
$\ringR= \fieldK[x,y,z]/(x^7+y^7-z^7)$, ${\rm char}\, \fieldK= p$.
\end{proposition}
\begin{proof}
For $p=2,23$ this was checked with the help of a computer and Remark
\ref{2,23remark}, so suppose that $p \neq 2,23$. Corollary
\ref{secondfrobenius27} then guarantees that $x^3 y^3 \notin
(x^4,y^4,z^4)^*$ if and only if $x^{3p^2} y^{3p^2} \notin (x^{4p^2},
y^{4p^2}, z^{4p^2})$. Write $p=7\ell + 2$ and $p^2=7k+4$ where
$k=7\ell^2 + 4\ell=p\ell + 2\ell$. Now $4p^2=7(4k+2)+2=28k+16$ and
so $z^{4p^2}$ equals $(x^7+y^7)^{4k+2} z^2$, so it is enough to show
that
$$x^{3p^2} y^{3p^2} \notin (x^{4p^2}, y^{4p^2}, (x^7+y^7)^{4k+2}) \, .$$
If this were not the case then we would have already, since
$\fieldK[x,y] \subset \ringR$ is a free extension,
$$x^{3p^2} y^{3p^2} \in (x^{4p^2}, y^{4p^2}, (x^7+y^7)^{4k+2}) \fieldK[x,y] \, ,$$
so we have to show that this is not true. By Lemma
\ref{lemmafrobenius1} we know that
\begin{equation*}
x^{3p} y^{3p} \notin (x^{4p}, y^{4p}, (x^7+y^7)^{4 \ell +1})
\end{equation*}
in $\fieldK[x,y]$. Since $\fieldK[x,y]$ is a regular ring, it is
$F$-pure, therefore we take a Frobenius power to conclude that
\begin{equation*}
x^{3p^2} y^{3p^2} \notin (x^{4p^2}, y^{4p^2}, (x^7+y^7)^{p(4 \ell
+1)}) \, .
\end{equation*}
But we have
$$p(4 \ell +1) =4 \ell (7 \ell +2) + 7 \ell +2 = 4(7 \ell
^2 +2 \ell) + (7 \ell +2) = 4(k - 2 \ell) + (7 \ell +2) = 4k -\ell
+2\, ,$$ which is strictly smaller than $ 4k +2$. Therefore
replacing the ideal generator $ (x^7+y^7)^{p(4 \ell +1)}$ by
$(x^7+y^7)^{4k+2} $ makes the ideal smaller, hence $x^{3p^2}y^{3p^2}
\notin (x^{4p^2}, y^{4p^2}, (x^7+y^7)^{4k+2})$ holds.
\end{proof}

\section{The case $ p\modeq 3 \modu 7$}

\begin{proposition}
\label{theorem 3 mod 7} If  $p\modeq 3 \modu 7$, then $x^3 y^3 \in
(x^4,y^4,z^4)^F \subseteq (x^4,y^4,z^4)^*$ in
$\fieldK[x,y,z]/(x^7+y^7-z^7)$ for ${\rm char}\, \fieldK =p$.
\end{proposition}
\begin{proof}
We show indeed that $x^{3p} y^{3p} \in (x^{4p}, y^{4p}, z^{4p})$.
Write $p=7\ell+3$; notice that $z^2 z^{4p}$ equals
$(x^7+y^7)^{4\ell+2}$, so it is enough to show that
$$x^{7(3\ell+1)+2} y^{7(3\ell+1)+2} = x^{3p} y^{3p} \in
(x^{4p}, y^{4p}, (x^7+y^7)^{4\ell+2}) = (x^{28\ell+12},
y^{28\ell+12}, (x^7+y^7)^{4\ell+2}) \, .$$ We will show that $
x^{7(3\ell+1)} y^{7(3\ell+1)} \in (x^{28\ell+12}, y^{28\ell+12},
(x^7+y^7)^{4\ell+2})$ holds in $\fieldK[x,y]$. Consider the
$(2\ell+1)\times(2\ell+1)$ matrix
$$A=\leftc {4\ell+2\choose{2\ell+1+i-j}} \rightc_
{1\leq i\leq 2\ell+1\atop{1\leq j\leq 2\ell+1}} \, .$$
Lemma \ref{lemma2} shows that
$$\det A=\prod_{t=0}^{2\ell}
\frac{{6\ell+2-t\choose{2\ell+1}}}{{2\ell+1+t\choose{2\ell+1}}} $$
and since $2\ell+1\leq 6\ell+2-t,2\ell+1+t < p$ for $0\leq t\leq
2\ell$ none of the binomial coefficients in the determinant vanishes
modulo $p$ and so $\det A$ is a unit modulo $p$. Now $A$, as a
matrix with entries in $\fieldK$, is invertible and we can find
$a_0, \dots a_{2\ell} \in \fieldK$ such that
$$A \left( \begin{array}{l} a_0\\ \vdots\\ a_{2\ell} \end{array} \right)
= \mathbf{e}_{\ell+1} \, ,$$
where $\mathbf{e}_{\ell+1}$ is the $\ell+1$th elementary vector of
size $2\ell+1$. Consider the polynomial
\begin{eqnarray*}
f&=& \left(\sum_{i=0}^{2\ell} a_i x^{7i} y^{7(2\ell-i)} \right) (x^7+y^7)^{4\ell+2}\\
&=& \left(\sum_{i=0}^{2\ell} a_i x^{7i} y^{7(2\ell-i)} \right)
\left(\sum_{j=0}^{4\ell+2} {{4\ell+2}\choose{j}}
x^{7j}y^{7(4\ell+2-j)} \right)\in \fieldK[x,y]
\end{eqnarray*}
and notice that the terms occurring in $f$ have the form $x^{7i}
y^{7(6\ell+2-i)}$ for $0\leq i\leq (6\ell+2)$. Working modulo
$x^{28\ell+12}, y^{28\ell+12}$, since
$$
\left\{
\begin{array}{lll}
7i& < &28\ell+12\\
7(6\ell+2-i)& < &28\ell+12
\end{array}
\right.
\Leftrightarrow\quad
2\ell+1 \leq i\leq 4\ell+1 ,
$$
we have
\begin{eqnarray*}
f &\equiv&  \sum_{i=2\ell+1}^{4\ell+1} \left( \sum_{j=0}^{2\ell} a_j {4\ell+2\choose{i-j}} \right) x^{7i}y^{7(6\ell+2-i)}\\
&=& \sum_{i=1}^{2\ell+1} \left( \sum_{j=1}^{2\ell+1} a_{j-1}
{4\ell+2\choose{2\ell+1+i-j}} \right) x^{7(i+2\ell)}y^{7(4\ell+2-i)}
\mod (x^{28\ell+12}, y^{28\ell+12}) .
\end{eqnarray*}
and our choice of $a_0, \dots, a_{2\ell}$ gives
$$ \sum_{i=1}^{2\ell+1}
\left( \sum_{j=1}^{2\ell+1} a_{j-1} {4\ell+2\choose{2\ell+1+i-j}}
\right) x^{7(i+2\ell)}y^{7(4\ell+2-i)}= x^{7(2\ell+\ell+1)}
y^{7(4\ell+2-\ell -1)}=x^{7(3\ell+1)} y^{7(3\ell+1)}
$$
and so $x^{7(3\ell+1)} y^{7(3\ell+1)} \in (x^{28\ell+12},
y^{28\ell+12}, (x^7+y^7)^{4\ell+2})$.
\end{proof}

\section{Conclusions and remarks}
\label{conclusions}

\noindent
Putting together the results of the previous sections we
obtain the following theorem.

\begin{theorem}
\label{mainTheorem}
Let $\fieldK$ denote a field of positive
characteristic $p$ and let $\ringR= \fieldK[ x,y,z]/(x^7+y^7-z^7)$.
Then $x^3y^3 \in (x^4,y^4,z^4)^*$ for infinitely many prime numbers
and $x^3y^3 \not\in (x^4,y^4,z^4)^*$ for infinitely many prime
numbers.
\end{theorem}
\begin{proof}
This follows directly from Propositions \ref{theorem 2 mod 7} and
\ref{theorem 3 mod 7}, taking into account Dirichlet theorem on
primes in an arithmetic progression, see for example \cite[Chapitre
VI, \S 4]{serrearithmetic}.
\end{proof}

\medskip
\noindent We can now settle the question posed by
M.~Hochster, C.~Huneke and the second author mentioned in the
introduction.

\begin{corollary}
There exists an ideal $J\subseteq \mathbb{Q}[x,y,z]/(x^7+y^7-z^7)$ which is tightly closed but
whose descents $J_p\subseteq \mathbb{Z}/p\mathbb{Z}[x,y,z]/(x^7+y^7-z^7)$ to characteristic $p$
are not tightly closed for infinitely many primes $p$.
\end{corollary}
\begin{proof}
Let $J$ be the tight closure in characteristic zero of the ideal
$(x^4,y^4,z^4)$ in $\mathbb{Q}[x,y,z]/(x^7+y^7-z^7)$; obviously
$J^*=J$. Since there are infinitely many primes $p$ satisfying
$p\modeq 2 \modu 7$, Proposition \ref{theorem 2 mod 7} shows that
$x^3 y^3 \notin J$. For the infinitely many primes $p$ satisfying
$p\modeq 3 \modu 7$ we have however $x^3 y^3\in (J_p)^*$ and so for
these primes $(J_p)^* \neq J_p$.
\end{proof}

\medskip\noindent
Surprisingly, we can also deduce from our considerations in positive
characteristic that the syzygy bundle $\Syz(x^4,y^4,z^4)$ is
semistable in characteristic zero (we do not know of a single prime
number where it is strongly semistable).

\begin{corollary}
\label{semistablechar0}
The syzygy bundle $\Syz(x^4,y^4,z^4)$ is
semistable on $C= \Proj \QQ[x,y,z]/(x^7+y^7-z^7)$.
\end{corollary}
\begin{proof}
Suppose that there exists a destabilizing sequence $0 \ra \shL \ra
\Syz(x^4,y^4,z^4)(6) \ra \shM \ra 0$, $\shL$ of positive and $\shM$
of negative degree. Such a sequence may be extended to a sequence on
the relative curve over an open subset of $\Spec \ZZ$. Let $c=
\delta(x^3y^3) \in H^1(C, \Syz (x^4,y^4,z^4)(6))$ denote the
cohomology class corresponding to $x^3y^3$ and let $c'$ denote the
image of $c$ in $H^1(C,\shM)$. If $c' \neq 0$, then its torsor would
be affine and $x^3y^3$ would not belong to the solid closure of
$(x^4,y^4,z^4)$ in characteristic zero. But then it would not belong
to the tight closure for almost all prime numbers (since affineness
is an open property), which contradicts Proposition \ref{theorem 3
mod 7}. Hence $c'=0$ and $c$ stems from a class $c'' \in H^1(C,
\shL)$. Modulo $p$, $c''_p$ is annihilated by a Frobenius power,
since $\shL$ has positive degree. But that would mean that also
$c_p$ would be annihilated by a Frobenius power and hence $x^3y^3
\in (x^4,y^4,z^4)^*$ for almost all prime numbers, which contradicts
Proposition \ref{theorem 2 mod 7}.
\end{proof}

\begin{remark}
How does $J= (x^4,y^4,z^4)^*$ in characteristic zero look like?
Since $\Syz(x^4,y^4,z^4)$ is semistable in characteristic zero by
Corollary \ref{semistablechar0}, we know that $\Syz(x^4,y^4,z^4)(m)$
is an ample sheaf for $m \geq 7$ and that the dual of
$\Syz(x^4,y^4,z^4)(m)$ is ample for $m \leq 5$. Since ampleness is
an open property, it follows for almost all prime numbers $p$ that
$R _{\geq 7} \subseteq (x^4,y^4,z^4)^*$ (even in the Frobenius
closure) and that $R_{\leq 5} \cap (x^4,y^4,z^4)^* \subseteq
(x^4,y^4,z^4)$. For degree $6$ we know that $x^3y^3, x^3z^3, y^3z^3
\not \in J$ by Proposition \ref{theorem 2 mod 7}. We do not know
whether $x^2y^2z^2$ and $xy^2z^3$ etc. belong to $J$ or not.
\end{remark}

\begin{remark}
What can we say in our example about tight closure and Frobenius
closure for the other remainders of $p$ modulo $7$? There is
numerical evidence showing that for $p \modeq 3,5,6 \modu 7$ the
element $x^3y^3$ belongs to the Frobenius closure of
$(x^4,y^4,z^4)$, but not for $p \modeq  1,2,4 \modu 7$. Moreover it
seems as if $x^3y^3 \in (x^4,y^4,z^4)^*$ for exactly $p \modeq
1,3,5,6 \modu 7$.

\smallskip\noindent
We began this work by looking at the example $xyz \in
(x^2,y^2,z^2)^*$ in $\fieldK[x,y,z]/(x^5+y^5-z^5)$. Here we have
strong computer evidence that $xyz \in (x^2,y^2,z^2)^F$ holds
exactly for the remainders $p \modeq 2,4 \modu 5$, and we have
proved this for $p \modeq 2 \modu 5$. Moreover, for $p \modeq 3
\modu 5$ we have proved as in Lemma \ref{notsemistable} and
Corollary \ref{secondfrobenius27} that the computation of tight
closure reduces to the question of whether $(xyz)^{p^2} \not\in
(x^{2p^2}, y^{2p^2}, z^{2p^2})$, but we were unable to settle this.
The difficulty lies in the fact that in reducing the statement to a
problem over $\fieldK[x,y]$ (and then to a matrix problem over
$\fieldK$), we have to replace $z$ twice, and have to deal with two
different kinds of binomial coefficients. For $p \modeq 1 \modu 5$
it is likely that $xyz \in (x^2,y^2,z^2)^*$ holds without being in
the Frobenius closure.
\end{remark}

\begin{remark}
It is known since the early days of tight closure that the Frobenius
closure $I^F$ of an ideal $I$ fluctuates arithmetically. The easiest
example is that $y^2 \in (x,y)^F$ holds in
$\fieldK[x,y,z]/(x^3+y^3+z^3)$ for prime characteristic ${\rm char
\,} \fieldK =p \modeq 2 \modu 3$, but not for $p \modeq 1 \modu 3$,
see \cite[Example 2.2]{hunekeparameter}. It is therefore not
surprising that our argument reduces the tight closure question to a
question about Frobenius closure.
\end{remark}

\begin{remark}
\label{solidclosure} Our example shows also that tight closure in
characteristic zero and in dimension two is not the same as solid
closure. Recall that an element $f$ in a local (or graded) excellent
domain $(\ringR, \fom)$ of dimension $d$ belongs to the solid
closure of an $\fom$-primary ideal $(f_1 \komdots f_n)$ if and only
if $H^d_{\fom} (\ringR[u_1 \komdots u_n]/(u_1f_1 \plusdots
u_nf_n+f)) \neq 0$ (see \cite{hochstersolid} and
\cite{brennertightproj}). In positive characteristic, tight closure
and solid closure are the same, and solid closure contains always
tight closure. The containment of $x^3y^3$ inside the solid closure
of $(x^4,y^4,z^4)$ in $\fieldK[x,y,z]/(x^7+y^7-z^7)$ follows from
Proposition \ref{theorem 3 mod 7} or from the fact that the syzygy
bundle is semistable in characteristic zero.

\smallskip\noindent
The example provides also an example of a ring $\ringR_{\ZZ}
=\ZZ[x,y,z]/(x^7+y^7-z^7)$  and an $\ringR_{\ZZ}$-algebra $A=
\ringR_{\ZZ}[u,v,w]/(ux^4 +vy^4+wz^4 +x^3y^3)$ such that $H^2_{\fom
\ringR_{\fieldK}}(A _{\fieldK})$ is zero for infinitely many prime
fields $\fieldK=\ZZ/(p)$ and non-zero for infinitely many prime
fields. The ring $A$ together with the ideal $\foa=(x,y,z)A \subset
A$ gives an example where the cohomological dimension of the open
subset $D(\foa)$ varies between $0$ and $1$ with the characteristic.
Classical examples for the dependence on the prime characteristic of
the cohomological dimension were given in
\cite[Example3]{hartshornespeiser} (see also
\cite[Corollary2.2]{singhwaltherarithmetic}), but as far as we know
our example is the first where it varies between $0$ and $1$,
corresponding to $D(\foa)$ being affine or not.
\end{remark}

\begin{remark}
\label{cohodimprojective} Let $Y \subset X$ denote a divisor on a
smooth projective variety over $\Spec \ZZ$ and let $Y_p \subset X_p$
denote the specialisations for a prime number $p$. How do properties
of $Y_p$ vary with $p\,$? Our example gives a smooth irreducible
divisor $Y$ on a smooth projective three-dimensional variety $X$
such that the complement $X_p -Y_p$ is an affine variety for
infinitely many but not for almost all $p$. Indeed, let $C= \Proj \,
\ZZ[x,y,z]/(x^7+y^7-z^7) \ra \Spec \ZZ$ be the relative curve and
let
$$0 \lra \shS= \Syz(x^4,y^4,z^4)(6) \lra \shS' =
\Syz(x^4,y^4,z^4,x^3y^3)(6) \lra \O_C \lra 0 $$ denote the extension
on $C$ defined by $x^3y^3$. Then $Y= \PP(\shS^\dual) \subset
\PP((\shS')^{\dual}) =X$ is a projective subbundle of codimension
one inside a projective bundle over $C$ of fiber dimension two. Our
result says that $X_p -Y_p =\PP((\shS'_p)^{\dual})
-\PP(\shS_p^\dual)$ is affine for $p \modeq 2 \mod 7$ and not affine
for $p \modeq 3 \modu 7$. We do not know whether such an example
exists if $X$ is a surface.
\end{remark}

\begin{remark}
\label{hilbertkunz} Our example is also relevant to the study of
Hilbert-Kunz multiplicities. The Hilbert-Kunz multiplicity is an
invariant of an ideal $I$ (primary to a maximal ideal) in a ring
$\ringR$ of positive characteristic $p$, defined by $e_{HK}(I)=
\lim_{e \in \NN} \length (\ringR / I^{[p^{e}]})/ p^{2e} \, \in \RR$,
see \cite[Chapter 6]{hunekeapplication}. It is related to tight
closure by the fact that $f \in I^*$ holds if and only if
$e_{HK}(I)= e_{HK}((I,f))$. Set $I=(x^4,y^4,z^4)$ and
$I'=(x^4,y^4,z^4,x^3y^3)$ in $\ZZ[x,y,z]/(x^7+y^7-z^7)$. Our results
give $e_{HK}(I_p)=e_{HK}(I_p')$ for $p \modeq 3 \modu 7$ and
$e_{HK}(I_p) \neq e_{HK}(I_p')$ for $p \modeq 2 \modu 7$. In
particular, the Hilbert-Kunz multiplicity is not eventually constant
as $p \ra \infty$.

\smallskip\noindent
On the other hand, V. Trivedi has shown in
\cite{trivedihilbertkunzreduction} that in the two-dimensional
graded situation the limit $\lim_{p \mapsto \infty} e_{HK}(I_p)$
exists. Moreover, one can show that this limit is the Hilbert-Kunz
multiplicity in characteristic zero as defined in
\cite{brennerhilbertkunzcriterion}. In our example we have $\lim_{p
\mapsto \infty} e_{HK}(I_p) =\lim_{p \mapsto \infty} e_{HK}(I'_p)$,
because they coincide for infinitely many prime numbers. This
corresponds to the fact that $x^3y^3$ belongs to the solid closure
of $I$ in characteristic zero. This limit is in our example $84$ (
see \cite[Introduction]{brennerhilbertkunzfunction} for the formulas
to compute the Hilbert-Kunz multiplicity). Apart from that we only
know for $p \modeq 2 \modu 7$ that $e_{HK}(I_p) \geq 84 + 28/ p^2$;
we get here only an inequality because the instability of
$\Syz(x^{4p},y^{4p},z^{4p})$ might be even worse that the
instability detected in Lemma \ref{notsemistable}.
\end{remark}

\begin{remark}
\label{nonstandard} H. Schoutens defined another variant of tight
closure for finitely generated algebras $R$ over $\CC$, called
non-standard tight closure (see \cite{schoutensnonstandard}). He
uses methods from model theory and an identification $\ulim
\overline{\ZZ/(p)} \cong \CC$, where $\ulim$ denotes the
ultraproduct with respect to a fixed non-principal ultrafilter. Then
the ultraproduct of the Frobenii of the approximations $R_p$ give a
characteristic zero Frobenius $R \ra R_\infty = \ulim R_p$ and yield
a new closure operation with several variants. A natural question is
whether these closure operations are independent of the choice of
the ultrafilter and whether the several variants coincide or not
(Question 1 after Theorem 10.4 in \cite{schoutensnonstandard}). Our
example shows at once that the so-called generic tight closure
depends on the choice of the ultrafilter. Moreover, if the parameter
theorem of Hara \cite[Theorem 6.1]{hunekeparameter} holds for
non-standard tight closure for two-dimensional graded $\CC$-domains,
then it follows that also non-standard tight closure depends on the
ultrafilter.
\end{remark}

\begin{question}
Suppose that $R$ is a finitely generated extension of $\ZZ$, let $I
\subseteq R$ denote an ideal and let $f \in R$. Set $M= \{ p \,
\mbox{ prime}:\, f_p \in (I_p)^* \}$. Is it possible to characterize
the subsets of the prime numbers which arise in this way? Do there
always exist congruence conditions which describe such an $M$ up to
finitely many exceptions?
\end{question}

\bibliographystyle{plain}

\bibliography{bibliothek}
\end{document}